\documentclass[11pt,reqno]{amsart}

\date{\today}

\usepackage[final]{graphicx}
\usepackage{amsfonts}
\usepackage{amsmath}
\usepackage{amssymb}
\usepackage{amsthm}
\usepackage{mathrsfs}
\usepackage[colorlinks=true,linkcolor=blue,citecolor=blue]{hyperref}
\usepackage{bbm}
\usepackage{setspace}

\newtheorem{theorem}{Theorem}

\newtheorem{corollary}[theorem]{Corollary}
\newtheorem{proposition}[theorem]{Proposition}

\newtheorem{definition}[theorem]{Definition}

\makeatletter\@namedef{subjclassname@2020}{\textup{2020} Mathematics Subject Classification}
\newcommand{\msc}[1]{\href{https://zbmath.org/classification/?q=cc:#1}{#1}}
\newcommand{\Email}[1]{\rm{\it E-mail:}\/~\href{mailto:#1}{\textsf{#1}}}

\newcommand{\R}{{\mathord{\mathbb R}}}

\newcommand{\Z}{{\mathord{\mathbb Z}}}

\renewcommand{\(}{\left(}
\renewcommand{\)}{\right)}
\newcommand{\irn}[1]{\int_{\R^n}{#1}\,dx}

\newcommand{\be}[1]{\begin{equation}\label{#1}}
\newcommand{\ee}{\end{equation}}

\newcommand{\Emph}[1]{\textit{\textbf{#1}}}

\begin{document}
\title[The importance of Luis Caffarelli's work in the study of fluids]{The importance of Luis Caffarelli's work in the study of fluids.}

\author[M.~J.~Esteban]{Maria J.~Esteban}
\address[M.~J.~Esteban]{CEREMADE (CNRS UMR No.~7534), PSL University, Universit\'e Paris-Dauphine, Place de Lattre de Tassigny, 75775 Paris 16, France. \newline\Email{esteban@ceremade.dauphine.fr}}


\begin{abstract}
In this paper we describe the work of Luis Caffarelli in the area of fluid mechanics and related topics. Not only has his work on fluid mechanics been very influential,
but many of his contributions that do not directly relate to fluid mechanics, such as his important results on the fractional Laplacian or the regularity of solutions to linear parabolic equations with oscillating coefficients, have been used in the study of fluids in many important ways. Thus, any review of his work has to include his contributions to the general (partial) regularity theory of solutions of Navier-Stokes equations and other studies related to fluid motion.

\vspace*{-0.5cm}
\end{abstract}



\thanks{\scriptsize \copyright~2024 by the author. Reproduction of this article by any means permitted for noncommercial purposes. \hbox{\href{https://creativecommons.org/licenses/by/4.0/legalcode}{CC-BY 4.0}}}

\maketitle
\thispagestyle{empty}



\smallskip
\section{Introduction}\label{Sec:intro}

Fluid mechanics is the branch of physics concerned with the mechanics of certain continuous media (liquids, gases, and plasmas) subjected to external forces. Its applications are quite broad: geophysics, oceanography, meteorology, astrophysics, biology, biophysics and aerospace, civil, chemical, and biomedical engineering, just to quote some of them.

\smallskip
Due to its practical interest, the study of fluid mechanics started very early. We can find discussions about it in ancient times: the treatise ``On floating bodies" written by Archimedes is considered to be the first study on fluid mechanics.  In Iran, the famous scientist Al Biruni was interested in experiments involving fluid motion; Leonardo da Vinci also made experiments and observations on fluids. But it is not until Newtonian mechanics was established that the study of fluids attained a more scientific and formal status. Since then, many mathematicians and mechanicians have been interested in fluid dynamics, like for instance Bernoulli, d'Alembert, Euler, Lagrange, Laplace, Poisson, etc, not to speak of all modern studies on mathematical fluids.

\smallskip
It has been known for a long time that pressure and viscosity are fundamental notions when  describing fluid flow. But it was only in the 19th century that the relationship between them was put into mathematical terms for incompressible and viscous fluids. The Navier-Stokes equations, named after Claude-Louis Navier and Gabriel Stokes, give an exact relationship between the viscosity of a fluid, the pressure, and the acceleration the fluid experiences. If you know the initial conditions of the flow and the external forces acting on it, then a solution of the equations will tell you how the fluid flow will behave at  future times.

\smallskip
Although it is not known how to solve the Navier-Stokes equations, one can show the existence of what is known as \it weak solutions\rm.  These are solutions that solve the equations, but in a very particular, weak, sense. They were introduced by Jean Leray in 1934. If one could show that those solutions, or at least some of them, were \it regular\rm, then the so-called Navier-Stokes problem would be solved. This is the central issue for the Navier-Stokes equations and a famous result of Caffarelli, Kohn and Nirenberg is the closest anybody has come to the study of the (partial) regularity of the \it weak solutions. \rm They showed that  \it weak solutions \rm enjoying some natural energy inequality inequality  are \it regular \rm except for a set of possible singularities that has $0$ one-dimensional dimension.

\smallskip
Luis Caffarelli  has made decisive contributions to mathematical analysis, in particular to the field of partial differential equations.  Probably his most important works are the study of fully nonlinear problems, local and nonlocal, the regularity of solutions of nonlinear partial differential equations and the analysis of the regularity of free boundaries. Likewise, his contributions on fluid dynamics, albeit  more qualitative than quantitative, have been very important. His work on regularity properties for fluids, both concerning the whole fluid domain, like the partial regularity results he proved with R. Kohn and L. Nirenberg \cite{MR0673830} and those concerning surface quasi-geostrophic flows he proved with A. Vasseur \cite{MR2680400}, the regularity of free boundaries in multi-phase problems involving fluids, his work on the Stefan problem, on jets and cavities with H.W. Alt and A. Friedman and on capillary drops with A. Mellet, are still considered the high points in this field.

\smallskip
The main mathematical questions  regarding fluids are the existence of solutions  as well as their main properties, like smoothness, (partial) regularity, or turbulence and chaotic behavior. Besides the theoretical analysis, there is the practical use of models to study fluids, in mechanics, in biology, in ecology, in weather forecast, in medicine, ... or in movie animation. Since most problems are too complex to be solved in a concrete and explicit manner, the numerical study of fluids has grown in importance, both concerning the numerical analysis of the codes used to produce approximate solutions of the equations and producing the approximate solutions themselves, often by very sophisticated numerical schemes.

\smallskip
As mentioned before, Caffarelli's work on fluids has been mainly related to the study of regularity of solutions and below we describe the main areas in which he has made seminal contributions. We say `seminal' because he (with collaborators) introduced a number of very novel and powerful methods that have been used later in many different contexts.
Here is a short timeline of his most important works in this area.

\smallskip
In his first big breakthrough, Caffarelli produced a proof for the regularity of melting ice \cite{MR0454350, MR0466965},  showing that the surface of a piece of ice remains mostly smooth as it turns into water. This result was only refined in 2024 by  A. Figalli and collaborators  \cite{MR4695505}.

\smallskip
Later,  Caffarelli got interested in models for fluids studying equations for their density, like for instance, in \cite{MR0560219}, where with A. Friedman he studied the properties of boundaries in rotating fluids in stars; or a little later his articles with H.W. Alt and A. Friedman about equations arising in jet flows and cavities for the stream function \cite{MR0637494, MR0647374, MR0818805}. The first paper contains an existence proof for the asymmetric two-dimensional jet problem for an ideal fluid with relaxed conditions on the nozzle shape. The next paper considers the same problem but with the added difficulty created by the presence of a gravity field. In both cases the shape of the flow was studied, but it was further described and studied in the third paper. In the three cases, the properties of the nozzle are of course of big importance. Other works with the same group of authors include \cite{MR0682265, MR0733897, MR0740956} and a comprehensive study of the compressive flows of jets and cavities can be found in  \cite{MR0772122}. With the same collaborators, Caffarelli also worked on  a related problem consisting in the study of the dam problem with two fluids \cite{MR0752593}. Again a free boundary separating the two fluids was considered.

\smallskip
More of less at the same time, Caffarelli co-authored a very important article about regularity in fluids, but this time the main subject of study was not a free boundary, but the full solutions of an equation modeling the dynamics of an incompressible viscous fluid, that is, the incompressible Navier-Stokes equation in $3$d. This was his very important article, with R. Kohn and L. Nirenberg, on the partial regularity of the Navier-Stokes equations \cite{MR0673830}. In a seminal work \cite{MR1555394} published in 1934, J. Leray demonstrated the existence of global weak solutions to the Navier-Stokes equations in three dimensions but the uniqueness and smoothness remain open. The partial regularity result of Caffarelli, Kohn and Nirenberg has had a huge impact in the literature of fluids, and it has basically not been improved so far.

\smallskip
Later on, Caffarelli considered fluids in other contexts, like for instance, studying capillary fluid drops with A. Mellet \cite{MR2307770, MR2373730}. The first article is concerned with the homogenization of a capillary equation for liquid drops lying on an inhomogeneous solid plane, while the second one deals with some properties of equilibrium liquid drops lying on a horizontal plane, when small periodic perturbations arise in the properties of that plane (due for example to chemical contamination or roughness).

\smallskip
The use of De Giorgi theory of regularity for solutions of nonlinear partial differential equations (see \cite{MR0093649, MR2680400}) allowed L. Caffarelli and A. F. Vasseur to prove the global regularity of the solutions to the critical surface quasi-geostrophic equations,  a simple model in climatology describing a fluid which is rapidly rotating. The approach they used to prove this regularity result was motivated by ideas of De Giorgi's
solution of the 19th Hilbert problem in 1957, which was devised to prove regularity for linear partial differential equations with measurable coefficients. The novelty of Caffarelli and Vasseur's approach was to apply it  in the context of strongly nonlocal operators. By doing so, they developed a theory which can be used in general situations like to study models of transport of scalars in turbulent fluid, phase transitions, or nonlocal image processing. Thus, not only the regularity result proved for the critical surface quasi-geostrophic equations was interesting, but also, importantly, the new method introduced to achieve it.

\smallskip
Later, with J.L. Vazquez, Caffarelli studied the porous medium flow with fractional potential pressure \cite{MR2847534}, and much later, again considering fractional diffusion, with M. Gualdani and N. Zamponi he proved the existence of weak solutions to some continuity equations with space-time nonlocal Darcy law \cite{MR4176917}.

\smallskip
What is striking in all the above works is that almost in each case the new and often surprising results were proved through the introduction of new methods, which would then be used for years and decades by many other mathematicians to study similar or different problems. Important examples of this are the remarkable \it monotonicity formula \rm that was used by Alt, Caffarelli and Friedman in \cite{MR0732100} to study the main properties of the solutions of two-phase free boundary problems. Another result of this kind is the introduction in \cite{MR0673830} of an iterative method producing the boundedness of the fluid velocity whenever the gradient is ``small" in a small neighborhood around that point. This method is called \it $\epsilon$-regularity. \rm We will give more details about this later.  The main point here is that the methods introduced in those papers have been used in different ways in the literature, but the basic idea behind them had their gestation in Caffarelli's and collaborators' work.

\smallskip
There are also aspects of the methods used in the above papers that have had a life of their own later on. For instance, in \cite{MR0673830} there was a need to control quantities related to the fluid velocity in weighted spaces. For that purpose the authors used a new set of interpolation inequalities, with weights, called nowadays \it Caffarelli-Kohn-Nirenberg inequalities\rm, which have been studied and used in many works and which are still in development.

\smallskip
Other problems related to fluids have been treated by Caffarelli in some of his papers, but here we shall restrain ourselves to dealing with the above topics.


\smallskip
To end describing the impact of Caffarelli's research on fluid mechanics, let us mention again that while some of his works were not about fluid mechanics, the methods he introduced in those articles have later influenced other people working on fluid mechanics. Let us describe two important examples of this. The first one is the work in collaboration with Silvestre on fractional Laplacians \cite{MR2354493}. The second one is the work in collaboration with Peral on $W^{1,p}$ estimates for elliptic equations in divergence form \cite{MR1486629}.
In the latter, Caffarelli and Peral introduced a quantitative approach to prove regularity of the solutions to equations with very oscillating coefficients. This method has had a lot of success in many applications in mechanics and in particular, it has been used to analyze the regularity of the stationary Navier-Stokes equations in domains with very low regularity properties. On the other hand, the fundamental work for fractional diffusion that is contained in \cite{MR2354493} has had for instance, applications in the study of the quasi-geostrophic equations.

\smallskip
The article is organized as follows. Section \ref{sec:NS} contains a description of Caffarelli, Kohn and Nirenberg's work about partial regularity for the solutions of Navier-Stokes equations in dimension $3$ as well as later related works. Section \ref{sec:QG} is devoted to the work of Caffarelli and Vasseur about global regularity for the critical quasi-geostrophic equations. In Section \ref{sec:jets} we describe the work of Caffarelli \it et al \rm on jets, cavities and other fluid  interface  problems.  The work of Caffarelli and collaborators on the regularity for the Stefan problem is the subject of Section \ref{sec:Stefan}. Section \ref{sec:drops} deals with results about the existence and regularity for capillary drop models. Finally, in Section \ref{sec:CKN} are described a number of results concerning interpolation inequalities with weights: the Caffarelli-Kohn-Nirenberg inequalities and extensions that were introduced in \cite{MR0673830} (see also \cite{MR0768824}).

\smallskip
\section{Partial regularity for the Navier-Stokes equations}\label{sec:NS}

\smallskip

The most common model to describe the dynamics of an incompressible fluid in a domain $\Omega$ of $\R^3$ are the Navier-Stokes equations for incompressible fluids with viscosity $\nu$:
$$u_t + u\cdot \nabla u - \nu\Delta u +\nabla p = f \quad\mbox{in}\quad (0, T)\times\Omega\,,$$
\vspace{-6mm}$$ \hspace{-60mm}\mbox{(NS)} \hspace{50mm} \nabla\cdot u= 0 \quad\mbox{in}\quad (0, T)\times\Omega\,,$$
\vspace{-5mm}$$ u(0,x)=u_0(x) \quad\mbox{in}\quad \Omega\,,$$
where $u(t,x):(0, T)\times \Omega \rightarrow \R^3 \,$ is the speed of the fluid at time $t$ and at the point $x\in \R^3$, $p$ stands for the pressure, $\nu$ denotes the fluid viscosity, $u_0$ is the initial velocity at time $t=0$,  and $f$ represents the sum of  the total external forces acting on the fluid. If $\Omega$ has a boundary, boundary conditions for $u$ on $\partial\Omega$ should also be made precise.

\smallskip
Note that the initial velocity $u_0$ has to satisfy conditions which are compatible with those imposed on $u$, that is, it should satisfy the same boundary conditions on $\partial\Omega$ and moreover, 
$$ \nabla \cdot u_0 = 0\quad\mbox{in}\quad \Omega\,.$$

\smallskip
It is very important to explain the meaning given to the above problem. If the functions $u$ and $p$ are smooth enough, $C^2$-regular for $u$ and $C^1$-regular for $p$, then the meaning of the problem is clear for all $x\in \Omega$. If the second derivatives of $u$ and the first derivatives of $p$ are defined almost everywhere, the equations will be understood in the sense of \it for almost every $x\in \Omega$. \rm But what if either $u$ or $p$ have less regularity?  Almost the weakest version that one can consider in this case is the notion of $(u,p)$ being a \it weak solution \rm of the above problem, that is, a solution in the sense of distributions: assume that $u_0\in L^2(\Omega)\,$ and $\, \nabla \cdot u_0 = 0\;\mbox{in}\; \Omega\,.$ It is said that $u\in L^\infty(0,T; L^2(\Omega))\cap L^2(0,T;H_0^1(\Omega)) $ is a \it weak solution \rm of (NS) if for all functions $\varphi\in C^\infty((0, T)\times\Omega, \R^3)$ with compact support in $(0,T)\times\Omega$, and satisfying $\nabla\cdot\varphi=0 \;\mbox{in } (0, T)\times\Omega$,
$$ \int_{(0, T)\times\Omega} \!\!\!\!\!\(- u_i\, \partial_t\varphi_i -u_i\,u_j \,\partial_{x_j}\varphi_i +\nu\, \nabla u_i \,\nabla \varphi_i\)dx\,dt = -\int_{(0, T)\times\Omega} f_i\,\varphi_i\,dx\,dt,\quad i=1,2, 3\,.$$
$$ \lim_{t\rightarrow0^+}\|u(t,\cdot) - u_0\|_{L^2(\Omega)}=0\,.$$

\smallskip
The existence of weak solutions to the Navier-Stokes equations was proved by Leray and Hopf in  \cite{MR1555394, MR0050423}. Leray dealt with the case $\Omega=\R^3$, while Hopf extended the result to the case of general domains $\Omega$. Moreover, Leray showed that these weak solutions are strong locally-in-time except on a certain compact set of singular times of Lebesgue measure zero.

\smallskip
Note that it is possible to prove that the weak solutions of (NS) satisfy the energy inequality
\be{energy-ineq} \int_{\{t\}\times \Omega}|u|^2\,dx + 2\, \int_0^t\int_\Omega |\nabla u|^2\,dx\,dt\le \int_\Omega |u_0|^2\,dx+ 2\int_0^t\int_\Omega f\cdot u\,dx\,dt\qquad \mbox{a.e. in} \,\,t\,,
\ee
 which are naturally satisfied by the smooth solutions of (NS).

\smallskip
If we want to speak about regularity of solutions $u$ of (NS) in more precise terms, let us say that a point $(x,t)$ is \it singular \rm if $u$ is not in $L^\infty_{\rm loc}$ in any neighborhood of $(x,t)$. The points that are not singular are called \it regular \rm  points, and they are the points where  $u$ is locally essentially bounded. What is called a  \it partial regularity theorem, \rm consists in proving an estimate for the dimension of the set $S$ of singular points.

\smallskip
After Leray's work, a long time passed without anybody managing to measure the size of the singular set. The first works in this direction are due to Scheffer (see \cite{MR0452123, MR0454426, MR0510154, MR0501249, MR0573611}). Scheffer's main result states the following
\begin{theorem}[\cite{MR0452123, MR0454426, MR0510154, MR0501249, MR0573611}]\label{thm:scheffer} 
 For $f=0$, there exists a  weak solution of the Navier-Stokes equations (NS) satisfying the local energy inequality whose singular set $S$ satisfies
$$\mathcal H^{5/3}(S) <\infty\,,$$
$$\mathcal H^1(S\cap(\Omega\times\{t\}))<\infty\quad\mbox{uniformly in}\; t,$$
where $\mathcal H^k$ denotes the Hausdorff $k$-dimensional measure, 
and where we say that a weak solution to (NS), $u$, satisfies the  \it local energy inequality  \rm if for almost every $t\in (0,T)$ and  for all $ \phi\in \mathcal D((0,T)\times\Omega)$, $\phi\ge 0$,

\be{local-energy-ineq}
\int_{\{t\}\times\Omega}\!\!\!\! | u|^2\phi\,dx+2\int_{(0,t)\times\Omega}\!\!\!\! |\nabla u|^2\phi\,dx\,dt\le \int_{(0,t)\times\Omega}\!\!\big(|u|^2(\phi_t+\Delta\phi)+(|u|^2 +2\,p)\,u\cdot\nabla\phi+2\,(u\cdot f)\,\phi\big)\,dx\,dt.
\ee
\end{theorem}

\smallskip
Scheffer's proof is inspired by the works of Mandelbrot \cite{MR0495674} and by the techniques used by Almgren to prove  partial regularity for minimal surfaces \cite{MR0420406}.

\smallskip
In a very celebrated paper \cite{MR0673830}, Caffarelli, Kohn and Nirenberg, improved Theorem \ref{thm:scheffer} proving a partial regularity result which basically remains up to now the best existing regularity result for (NS). In order to state it, let us first define a notion of weak solution that they introduced in this paper.

\begin{definition}\label{def:suitableCKN}
\rm A weak solution to (NS) in an open domain $D$ of $\R_+\times\R^3$, $(u,p)$, is called \it suitable \rm  if  $u, p$ and $f$ are measurable functions in $D$,  $p\in L^{5/4}(D)$, $f\in L^q(D)$ for some $q>5/2$, $\nabla\cdot f=0$ in $D$, if $u$ has finite energy,  that is,
$$ \iint_D |\nabla u|^2dx\,dt<\infty\,;$$
if there exists $E_0>0$ such that $\int_{D_t} |u|^2dx\le E_0$,  $D_t:=D\cap\(\{t\}\times\R^3\)$, 
for almost every $t$ such that $D_t\neq \emptyset$ and 
if the local energy inequality \eqref{local-energy-ineq} is satisfied.
\end{definition}
Their main theorem states the following.
\begin{theorem}[\cite{MR0673830}]\label{thm:CKN-NS}
Let $\Omega$ be an open set of $\,\R^3$, and  $u_0\in H^1_0(\Omega)$ satisfy
$$ u_0=0\;\mbox{on}\;\; \partial\Omega\,,\;  \;\; \nabla\cdot u_0=0\;\mbox{ in}\;\; \Omega$$
in a suitable weak sense. Let $f\in L^2((0,T)\times\Omega)$. If $u$ is a \it suitable weak solution \rm to (NS) with
$$\,u(t, \cdot)=0\;\mbox{on}\;\; \partial\Omega\,,\; \forall t\in (0, T],$$
then $u$ is bounded outside a singular set $S$, such that
$$ \mathcal P^1(S)=0\,,$$
where $\mathcal P^1$ is a measure  on $\R\times\R^3$ analogous to the one-dimensional Hausdorff measure, but defined using parabolic cylinders instead of Euclidean balls (in particular, $\mathcal H^1 \le C\,\mathcal P^1$ for some positive constant $C$).
\end{theorem}
Theorem \ref{thm:CKN-NS} is a big improvement of Theorem \ref{thm:scheffer} because  it allows external forces which are not equal to $0$ and it improves on the Hausdorff measure for which an estimate of the singular set of $u$ is secured. And also, this theorem has a \it local \rm character.

\smallskip
The techniques introduced to prove Theorem \ref{thm:CKN-NS} have had a big impact, having been used once and again in many later articles, by practically all the specialists on regularity for the evolution Navier-Stokes equations.

\smallskip
Let us go through the main arguments in the proof of Theorem \ref{thm:CKN-NS}. The first very important remark which the paper contains is the \it scalability \rm of (NS), that is, the fact that if $u$ and $p$ are solutions to (NS),  for each $\lambda>0$,
$$ u_\lambda(t,x):= \lambda\,u(\lambda^2 t, \lambda\,x)\,,\quad  p_\lambda(t,x):= \lambda^2\,p(\lambda^2 t, \lambda\,x)$$
are solutions to (NS) with $f$ replaced by $ f_\lambda(t,x):= \lambda^3\,f(\lambda^2 t, \lambda\,x)$. Since the above shows that in the framework of (NS) 'time' $t$ has dimension $2$, while the space variables $x_i $ have dimension $1$, it is natural to use parabolic cylinders
$$ Q_r(t,x):= \{(\tau, y)\,:\; t-r^2<\tau<t+r^2,\; |y-x|< r\}$$
as natural neighborhoods of any space-time point $(t,x)$.

\smallskip
As the authors explain in their article, ``it is well known that if $u_0$ and $f$ are small enough in a suitable norm, then the solution of (NS) stays regular at least for a short time. The main task in proving a partial regularity theorem lies in proving a local, dimensionless version of an analogous result". They do it in two steps. The first consists in proving a local version of the main result in \cite{MR0510154}. It states that  if $u, p$, and $f$ are ``sufficiently small" on the unit cylinder $Q_1:=  Q_1(0,0)$, then $u$ is regular on the smaller cylinder $Q_{1/2} := Q_{1/2}(0,0)$:
\begin{proposition}[\cite{MR0673830}]\label{prop:P1} There are absolute constants $\varepsilon_1$ and $C_1>0$, and a constant $\varepsilon_2(q)> 0$ depending only on $q$, with the following property. Suppose $(u,p )$ is a \it
suitable weak solution \rm of the Navier-Stokes system on $Q_1$ 
with force $f\in L^q$, for some $q>5/2$; suppose further that
$$M(1):=\int_{Q_1}\(|u|^3 + |u||p|\)\,dx\,dt + \int_{-1}^0\(\int_{|x|<1}|p|\,dx\)^{5/4}dt\le \varepsilon_1\,,$$
and
$$F_q(1):= \int_{Q_1}|f|^q\,dx\,dt\le \varepsilon_2\,.$$
Then
$$|u(t,x)|\le C_1$$
for Lebesgue-almost every $(t,x)\in Q_{1/2}$. 
In particular, $u$ is regular in $Q_{1/2}$.
\end{proposition}

\smallskip
The above result is a local regularity result, but does not contain any information about the size of the singular set of $u$. But one could use it, as done by Scheffer, to prove Theorem \ref{thm:scheffer}.  What Caffarelli, Kohn and Nirenberg did instead was to use the scalability property of (NS) to consider $M(r)$ and $F_q(r)$, $r>0$,  instead of $M(1)$ and $F_q(1)$ to prove the following:
\begin{corollary}[\cite{MR0673830}]\label{cor:C1}
Suppose that $(u,p)$ is a \it suitable weak solution \rm of (NS) on some cylinder $Q_r$ with force $f\in L^q$, $q>5/2$. If moreover
$$ M(r) \le \varepsilon_1\,,\;\; F_q(r)\le \varepsilon_2\,,$$
then
$$|u(t,x)|\le C_1\,r^{-1}\,,$$
Lebesgue-almost everywhere in $Q_{r/2}(t,x)$. 
\end{corollary}
Then, they used this rescaled version of Proposition \ref{prop:P1} to prove the second main result leading to the partial regularity result stated in Theorem \ref{thm:CKN-NS}:
\begin{proposition}[\cite{MR0673830}]\label{propo:P2}
There is an absolute constant $\varepsilon_3$ such that if $(u,p)$ is a \it suitable weak solution \rm of (NS) near $(t,x)$ and if
$$ \limsup_{r\rightarrow 0^+}\;\frac1r\int_{Q_r(t,x)}|\nabla u|^2\,dx\,dt \le \varepsilon_3\,,$$
then $(t,x)$ is a regular point of $u$.
\end{proposition}
A standard covering argument allows then to prove Theorem \ref{thm:CKN-NS}.

\smallskip
Note that the above proposition is what is called an \it $\varepsilon$-regularity \rm theorem. In later works about partial regularity results, other $\varepsilon$-regularity \rm theorems have been proved or used: for instance, in \cite{MR1488514}, Lin  proved  \it $\varepsilon$-regularity \rm choosing a smallness condition on 
\be{eps-Lin}
\int_{Q_1}\(|u|^3 + |p|^{3/2}\)\,dx\,dt\,,
\ee
 while in \cite{MR2548880} Wolf choose to prove partial regularity using $\varepsilon$-regularity applied to the quantity
\be{eps-Wolf}\limsup_{r\rightarrow 0^+}\;\frac1r \int_{Q_r(t_0, x_0)}\Big|\mbox{curl } u \times \frac{u}{|u|}\Big|^2\,dx\,dt\,.
\ee


\medskip
The main novel features concerning partial regularity of the solutions to (NS) to be found in the article of Caffarelli, Kohn and Nirenberg \cite{MR0673830} are the following: 

\smallskip
-- first, a very explicit definition of the notion of \it suitable weak solutions\rm, first introduced by Scheffer, and which has been used once again since then in the works of many specialists, and remains still important, because it is  the good notion, the most practical one, for the study the local regularity of solutions to Navier-Stokes equations. Note that in their proof Caffarelli, Kohn and Nirenberg did not use the (NS) equations directly. They mainly used  the local energy inequalities and a formula for the pressure that allowed them to get the final partial regularity result.

\smallskip
-- Then, the idea to choose the smallness of some \it critical quantities \rm in a neighborhood of a given Lebesgue point as a sufficient condition for proving the boundedness of the solutions at that point, and therefore the local regularity. This was somehow already present in Scheffer's papers, but was pushed much further in \cite{MR0673830}. The locality of the procedure was a key point in Caffarelli-Kohn-Nirenberg's work.

\smallskip
-- Finally, the proof of the partial regularity of \it suitable weak solutions \rm to the Navier-Stokes equations, result that has remained unchallenged until now, which is not surprising because this result is optimal for the \it suitable solutions\rm, that is, for the solutions satisfying the local energy inequality. And up to now, nobody has produced another regularity that goes beyond this for general solutions, unless one makes other assumptions on them. 

\smallskip
Caffarelli, Kohn and Nirenberg figured out how to use dissipation to reduce the dimension of the singular set with respect to Scheffer's result.  Indeed, they proved for the first time that if locally, the integral of the square of the velocity gradient is small in a little ball around some Lebesgue point of $u$, then there is local regularity near the center of that ball. This is not trivial at all, because  controlling the velocity gradient at a given point does not \it a priori \rm control how big the velocity could be. They figured out a very clever way to use scaling to deal with that problem. \it A priori \rm one could have a huge value of $u$ near that point, but going to a smaller ball,  that value is decreased. And  this procedure can be iterated again and again until getting a good estimate for $u$. Then one has to average all the information gotten in this iteration to see that control has not been lost in the iterative process.  This procedure is certainly the major contribution of the paper, and it has been used in practically all later partial regularity results. Different quantities can be used  to control the iteration, but the iteration will be needed. One can organize the proof differently and even make it more direct and easy to follow, but the main ingredient will still be the iteration on the radiuses of the ball until one manages to control the velocity in a small neighborhood of a regular point.

\smallskip
The $\varepsilon$-regularity argument is the main reason why the work of Caffarelli, Kohn and Nirenberg published in \cite{MR0673830} remains to date an unavoidable reference when working on partial regularity for the Navier-Stokes equations. And until something radically new is discovered, the result in \cite{MR0673830} will not be challenged because it is optimal in the framework of solutions satisfying the local energy inequality.

\smallskip
There have been other proofs of very close versions of the partial regularity theorem proved in \cite{MR0673830}. First, there was the work of Lin in \cite{MR1488514}, still for $\Omega=\R^3$ and for $f\equiv 0$. His paper simplifies the proof in \cite{MR0673830} by choosing \eqref{eps-Lin} as critical quantity for the $\varepsilon$-regularity result. Therefore, the partial regularity result proved by Lin in \cite{MR1488514} makes a stronger assumption on the pressure, $p\in L^{3/2}$,  than what was assumed in \cite{MR0673830}. He could do it thanks to the regularity for the pressure proven by Sohr and von Wahl \cite{MR0847086} in 1986, four years after the Caffarelli-Kohn-Nirenberg's paper was published.  Motivated by the choice of a different quantity to measure the smallness of a solution in a neighborhood of some point, Lin modified the definition of \it suitable weak solution \rm  by replacing in it the condition $p\in L^{5/2}$ by $p\in L^{3/2}$. Then, he proved a compactness result for bounded \it suitable weak solutions \rm which allowed him to simplify the proof of Caffarelli, Kohn and Nirenberg.

\smallskip
In 2007  Vasseur published another quite different proof of the partial regularity of solutions to the incompressible Navier-Stokes equation in dimension 3 (see \cite{MR2374209}) using a method introduced by De Giorgi for proving regularity of solutions to linear elliptic equations with rough diffusion coefficients. See also Caffarelli-Vasseur-2010 \cite{MR2660718}.

\smallskip
Some years later, the paper by Wolf \cite{MR2548880} used other  other critical quantities, \eqref{eps-Wolf}, that is, a different $\varepsilon$-regularity result,  to somehow simplify the proof of Caffarelli, Kohn and Nirenberg.

\smallskip
An important change of direction came with the paper by Jia and Sverak \cite{MR3056752} who studied the solutions of (NS) perturbed by a subcritical drift. They showed an $\varepsilon$-regularity result for that problem that was key in the proof of existence of self-similar solutions with big initial data, potentially leading to the proof of non-uniqueness for the Leray solutions of the problem.  This important extension of the work contained in \cite{MR0673830} has been a source of inspiration for many researchers. For instance,  in \cite{MR4076070} Barker and Prange considered recently the case of (NS) equations perturbed by a critical drift and also showed an $\varepsilon$-regularity result that lead to a local regularization result in space and time. More precisely, if one starts with an initial datum in $L^2$ which is locally in $L^3$, can one prove a local regularity result in that region in space and time. Inspired by the results in \cite{MR3056752}, more or less at the same time, Kang-Miura and Tsai proved in \cite{MR4379146} a similar result using instead a compactness approach \it \`a la \rm Lin.

\smallskip
Another recent interesting development appeared in the work of Albritton, Barker and Prange: in \cite{MR4593277} there is a new concise proof of a certain one-scale $\varepsilon$-regularity criterion using weak-strong uniqueness for solutions of the Navier-Stokes equations with non-zero boundary conditions, inspired by an analogous approach for the stationary system due to Struwe (see \cite{MR0933230}). In this paper, Albritton, Barker and Prange use an $\varepsilon$-regularity result based on the quantity $|u|^4$, that is, they look for smallness of 
$$
\int_{(0,T)\times \Omega}|u|^4\,dx\,dt\,,
$$
and in this case they prove both a qualitative and a quantitative statement of $\varepsilon$-regularity.

\smallskip 
Theorem \ref{thm:CKN-NS} is not the only important result contained in  \cite{MR0673830}. Indeed, apart from the estimate of how big the singular set can be, the local theory of regularity developed in the paper allows to study the initial value problem in $\R^3$ for $f=0$, with information about the possible location of the singular set. They state the following.
\begin{theorem}[\cite{MR0673830}] \label{thm:THEOREMC} Suppose $u_0\in L^2(\R^3)$ satisfies $\nabla\cdot u_0= 0$, and suppose 
\be{cond1} \int_{\R^3} |x||u_0|^2\,dx<+\infty.\ee
Then there exists a weak solution of the initial value problem (NS) with $f = 0$ which is regular in the region
$$\{ (t,x)\,:\; |x|^2t > K_1\},$$
where $K_1$ is a constant depending only on $ \int_{\R^3} |x||u_0|^2\,dx$ and on  $\int_{\R^3} |u_0|^2\,dx$.
\end{theorem}

Condition \eqref{cond1} says that $u_0$ ``is not  too singular" at the origin. The next result in \cite{MR0673830} states that under suitable conditions,  
a \it suitable weak solution \rm $u$ with initial data satisfying  \eqref{cond1} has
the property 
\be{cond2}
\sup_{0<\tau<t}\int_{\{\tau\}\times\R^3} |u|^2|x|^{-1}< +\infty\,,\quad \int_0^t\int_{\R^3} |\nabla u|^2|x|^{-1}<+\infty \;\mbox{for each } t\,.
\ee
Using this, the next theorem shows that the line $\{x = 0\}$ consists entirely of regular points:  $u$ is regular in a parabola above the origin.
\begin{theorem}[\cite{MR0673830}]\label{thm:THEOREMD} There exists an absolute constant $L_0>0$ with the following property. If $u_0\in L^2(\R^3)$ with $\nabla\cdot u_0= 0$, and if \eqref{cond1} holds, and if
$L:=\int_{\R^3} |u_0|^2|x|^{-1}\,dx<L_0$, then there exists a weak solution of the initial value problem (NS) with $f = 0$ which is regular in the region
$$\{(t,x ) : |x|^2<t\,(L_0-L)\}\,.$$

\end{theorem}

According to the results leading to the proof of Theorem \ref{thm:CKN-NS}, and in particular Proposition \ref{propo:P2}, a \it suitable weak solution \rm to (NS), $u$,  is regular at $(t,x)$ if either it is locally ``sufficiently small" in a suitably scaled sense, or if $|\nabla u|^2$ is not too singular at $(t,x)$. Then, Caffarelli, Kohn and Nirenberg thought of using this information in conjunction with additional hypotheses on the initial data $u_0$, to get extra information on the set of possible singularities of the solution to the initial value problem. They did that in the case $\Omega=\R^3$ and $f=0$.

\smallskip
In order to prove the two above results, and because they needed to control quantities involving $u$ and $\nabla u$  in some weighted $L^p$ spaces,  Caffarelli, Kohn and Nirenberg introduced some interpolation inequalities with weights that have lived an independent life of their own, independently of the use Caffarelli, Kohn and Nirenberg had for them in their paper. The main reason why the weights appeared in these proofs, is because when dealing with the definition of \it suitable weak solutions, \rm \eqref{local-energy-ineq}, they chose test fonctions $\phi$ which behaved like $|x|$ or $|x|^{-1}$ in some regions of $\R^3$.

\smallskip
These inequalities read as follows:
\be{CKN-ineq-orig1}
\(\irn{|x|^{-b p} |u|^p}\)^{1/p}\leq C\(\irn{|x|^{-a r}|\nabla u|^r}\)^{c/r}\(\irn{|x|^{\beta q}|u|^q}\)^{(1-c)/q}\,,
\ee
with 
$$ q,r\ge 1\,,\;\;p>0\,,\;\; 0\le c\le 1\,,\;\; -b= c\,\sigma + (1-c)\beta\,,$$

$$\frac1r - \frac{a}{n}\,,\quad \frac1q+\frac\beta{n}\,,\quad \frac1p- \frac{b}{n} >0\,.$$
In the subsequent article \cite{MR0768824} the authors developed more in detail the theory of the above inequalities, and in particulier they described the sets of parameters $\alpha, \beta, \gamma, \sigma$ for which the constant $C$ in the above inequality is finite. We will delve on these inequalities more in detail in Section \ref{sec:CKN}. 

\smallskip
To finish with this long section on partial regularity results emanating from the original article of Caffarelli-Kohn and Nirenberg, let us note that something that was not done in this paper, and that has remained not understood until now, is the possible relation between the singularity set and turbulence, a very important phenomenon whose current understanding is far from begin satisfactory.

\smallskip
\section{Regularity for the solutions to the critical surface quasi-geostrophic equations}\label{sec:QG}

In oceanography and meteorology a fundamental equation in the context of a rapidly rotating, density stratified, viscous, incompressible fluid is the $3$-dimensional Navier-Stokes equation. Both the forces of rotation and stratification impose a tendency toward $2$-dimensionality on the $3$-dimensional fluid motion. Some simplifying assumptions reduce the problem to the study of a $2$-dimensional equation which describes the evolution of the temperature field on a surface that bounds the fluid. In the geophysical fluids literature this equation is known as the surface quasi-geostrophic equation.

\smallskip
The surface quasi-geostrophic active scalar equation (QG) (\cite{MR1267050, MR1304437, MR2196360}), without dissipation, is a model of geophysical origin which started to be studied mathematically mainly because of its similarity to the $3$d Euler equations and because, as a scalar model in $2$d, it is simpler to study compared to the full $3$d Euler equation. The equation for (QG) is:
\be{QG1}  \partial_t \theta + u \cdot \nabla \theta = 0\,,
\ee
with
\be{u-theta} u = R^\perp \theta\,,
\ee
where $R := \nabla \Lambda^{-1}$ is the Riesz transform, $\Lambda: = (-\Delta)^\frac12$ being the Zygmund operator.  The second problem is close to this one, but including dissipation. The corresponding equation is 
\be{QG2} \partial_t \theta + u \cdot \nabla \theta +\kappa\,\Lambda^s\theta= 0\,,
\ee
where $u$ and $\theta$ are linked as in \eqref{u-theta}, $\kappa>0$, $s>0$. When $s=1$, the problem is called \it critical \rm and it corresponds to friction with boundaries. If $s>1$ the case is called \it subcritical \rm and \it supercritical \rm when $s<1$. 

\smallskip
Note that the above equation involves fractional diffusion, which already appeared before  in problems related to  surface flame propagation and in financial mathematics. We will describe important work of Caffarelli \it et al \rm on fractional diffusion later on  in the paper.

\smallskip
As stated in \cite{MR3187680}, the dissipative SQG equation \eqref{QG2} describes the evolution of a surface temperature field $\theta$  in a rapidly rotating, stably stratified fluid with potential vorticity \cite{MR1304437, MR1312238}. From the mathematical point of view, the non-dissipative SQG equations (QG), that is \eqref{QG1}-\eqref{u-theta},  have properties that are similar to those of the $3$d Euler equations in vorticity form \cite{MR1304437} and yet one may for instance prove the global existence of finite energy weak solutions \cite{Resnick}, but for completely different reasons than for $2$d Euler. This feature enhanced greatly the interest around the study of the surface quasi-geostrophic equations.

\smallskip 
The existence of weak solutions to \eqref{QG2}-\eqref{u-theta} was proved in \cite{Resnick} and it was  known very early that the subcritical equations have
smooth solutions (\cite{MR2215601, MR1709781}). The critical and supercritical cases are much more difficult to deal with.  In the critical case, it was known, see \cite{MR1855665, MR2084005},  that if the initial data was small in $L^\infty$, then the solution would be global and regular for all times. But what about all other cases? 

\smallskip
In both the critical and supercritical cases, Chae and Lee considered in \cite{MR1962043} the well-posedness with initial conditions small in Besov spaces (see also Wu \cite{MR2111923}).

\smallskip
Then, around 2006-2007, two remarkable and quite different proofs of global existence and regularity for the critical case were published. One of them was contained in   \cite{MR2680400}, where Caffarelli and Vasseur used harmonic extension to prove a gain of regularity for Leray-Hopf weak solutions, in the spirit of De Giorgi's technique for uniformly elliptic equations with measurable coefficients. First they proved $L^\infty$ estimates starting from $L^2$ ones and then,  $C^\alpha$ regularity. In essence, they proved H\"older regularity for the solutions of \eqref{QG2} whenever the initial datum is in $L^2$. The other result was proved by Kiselev, Nazarov and Volberg in \cite{MR2276260}, article which contains a quite different and very elegant proof of the fact that in $2$d, solutions with periodic $C^\infty$ data for the critical surface quasi-geostrophic equations remain $C^\infty$ for all times. Their main argument was a non-local Maximum Principle for a suitably chosen modulus of continuity. In both cases, the external force $f$ was equal to $0$. Some years later, a third proof of global regularity in the critical case was published by Kiselev and Nazarov in \cite{MR2749211}. This new proof used an appropriate set of test functions and estimates on their evolution. 

\smallskip
Let us explain in some detail the result proved by Caffarelli and Vasseur  in \cite{MR2680400} and the method used in the proof.

\smallskip
The first result is concerned with a gain of regularity from $L^2$ to $L^\infty$:
\begin{theorem}[\cite{MR2680400}]\label{thm:thmQG1} Let $\theta\in L^\infty(0, T; L^2(\R^N))\cap L^2(0,T; H^{1/2}(\R^N))$ and for every $\lambda > 0$ define:
$$\theta_\lambda =(\theta-\lambda)_+\,.$$
If  for every $\lambda > 0$, $\pm\,\theta_\lambda$ satisfy the level set energy inequalities:
\be{level-set-energy-ineq}
\int_{\R^N} \theta_\lambda^2(t_2, x)\,dx +2\int_{t_1}^{t_2}\int_{\R^N}|\Lambda^{1/2}\theta_\lambda|^2\,dx\,dt \le \int_{\R^N} \theta_\lambda^2(t_1, x)\,dx \;,\; \forall \,0<t_1<t_2\,,
\ee
then, 
$$ 
\sup_{x\in \R^N}\theta_\lambda(T, x) \le C\,\,\frac{\|\theta_0 \|_{L^2(\R^N)}}{T^{N/2}}\,,
$$
where $C$ is a positive constant.
\end{theorem}
Note that the level set energy inequalities \eqref{level-set-energy-ineq} are naturally expected to hold for solutions of \eqref{QG2}-\eqref{u-theta}. Indeed, it is enough to write the operator $\Lambda$ as the normal derivative of the harmonic extension of $\theta$ to the upper half space. This idea was extended to treat general fractional derivative operators by Caffarelli and Silvestre in \cite{MR2354493}.

\smallskip
With some additional information about $v$ one can get extra regularity, from $L^\infty$ to $C^\alpha$:
\begin{theorem}[\cite{MR2680400}]\label{thm:thmQG2} 
For $r>0$,  define $Q_r:= [-r, 0]\times[-r, r]^N$. If a function $\theta(t,x)$ satisfying \eqref{level-set-energy-ineq} is bounded in $[-1, 0]\times\R^N$ and if $v_{| _{Q_1}}\in L^\infty(-1, 0, \mbox{\sl BMO}(\R^N))$, then $\theta\in C^\alpha (Q_{1/2})$.
\end{theorem}
Then the main result about solutions to \eqref{QG2}-\eqref{u-theta}  is proved:

\begin{theorem}[\cite{MR2680400}]\label{thm:thmQG3} Let $\theta$ a solution to \eqref{QG2}-\eqref{u-theta} satisfying the level set energy inequalities \eqref{level-set-energy-ineq}. Then, for every $t_0>0$,  there exists $\alpha\in (0,1)$ such that $\theta\in C^\alpha([t_0, \infty)\times \R^N)$.

\end{theorem}
The above result  is a straightforward corollary of the two above theorems. Indeed, Theorem \ref{thm:thmQG1} implies that $\theta$ is uniformly bounded on $[t_0, \infty)$ for every $t_0 > 0$. Since singular integral operators are bounded from $L^\infty$ to $ \mbox{\sl BMO}$,  $u \in  L^\infty(t_0, \infty;\, \mbox{\sl BMO}\,(\R^N))$ and, after proper scaling, Theorem \ref{thm:thmQG2} implies the result of Theorem \ref{thm:thmQG3}.

\smallskip
Furthermore, note that higher regularity for solutions to \eqref{QG2}-\eqref{u-theta} follows from standard potential theory.

\smallskip
\noindent\bf Remark. \rm A key ingredient in the proofs of Theorems 1 and  2 is De Giorgi's ``oscillation lemma" (see \cite{MR0093649}) to prove H\"older continuity. Indeed, as Caffarelli and Vasseur explain in their article, \it the techniques in \cite{MR2680400}
could be seen as a parabolic De Giorgi-Nash-Moser method to treat \it linear boundary parabolic
problems \rm of the type:
$$
\mbox{div}(a\nabla \theta)=0 \;\;\mbox{in }\;  \Omega\times [0,T]\,;\quad  [f(\theta)]t =\theta_\nu \;\mbox{ on }\; \partial\Omega\times[0,T],
$$
that arise in boundary control theory (see Duvaut-Lions \cite{MR0464857, MR0521262}). \rm

\smallskip
As already noted in other parts of this article,  using the same \it linear techniques, \rm similar results to Theorem \ref{thm:thmQG3} can be obtained even for systems, in particular in applications in fluid mechanics (see the new proof of Caffarelli-Kohn-Nirenberg's partial regularity result of Vasseur in  \cite{MR2374209} and the work on $L^p$ estimates for quantities advected by compressible flows of Mellet and Vasseur \cite{MR2521733}).

\smallskip
The critical case in bounded domains has been studied in several works, the first one being \cite{MR3595454}.

\smallskip
In the supercritical range $s < 1$ , the global regularity 
of the Leray-Hopf weak solutions to the surface quasi-geostrophic equations is an open problem related to the global existence of classical solutions: in fact, it is well-known that
Leray-Hopf weak solutions coincide with classical solutions as long as the latter exist.
Constantin and Wu \cite{MR2466323, MR2483817} obtained partial results by extending the program of Caffarelli and Vasseur in \cite{MR2680400} to the
supercritical regime. 

\smallskip
For some time, the only regularity results available for large initial data were on conditional regularity and finite time regularization of solutions. For instance, it was shown by Constantin and Wu in \cite{MR2466323} that if the solution is $C^\delta$ with $\delta > 1 - \alpha$, then it is smooth (Dong-Pavlovic \cite{MR2525640} later improved this result to $\delta = 1 -\alpha$). Finite time regularization was proved by Silvestre \cite{MR2595196} for $\alpha$ sufficiently close to $1$, and for the whole dissipation range $0 < \alpha < 1$ by Dabkowski \cite{MR2773101} (with an alternative proof of the latter result given in \cite{MR2803787}). Also,  there are several results proving \it eventual regularity, \rm that is, regularity for large times, see for instance \cite{MR2595196, MR2803787, MR2773101, MR3498176}.
In contrast, a recent interesting result by Buckmaster, Shkoller and Vicol in \cite{MR3987721} establishes the non-uniqueness of a class of (very) weak-solution for the system \eqref{QG2}-\eqref{u-theta}, even for subcritical dissipations (see also the more recent article  by Cheng, Kwon and Li \cite{MR4340931}). Other interesting loss-of-regularity theorems have been proved by C\'ordoba and Martinez-Zoroa: in \cite{MR4452676}  they prove the ill-posedness of the inviscid quasi-geostrophic equations  \eqref{QG1}-\eqref{u-theta} in the spaces $C^k$ and also in Sobolev spaces $H^s$ for $s\in (3/2, 2)$, while in  \cite{CMZ-2} they consider the supercritical dissipative surface quasi-geostrophic equations \eqref{QG2}-\eqref{u-theta} for which they prove the existence of a unique global solution  which loses regularity instantly.

\smallskip
Recently,  again in the supercritical case, Colombo and Haffter published in \cite{MR4235801} an estimate on the dimension of the singular set of the solutions to \eqref{QG2}-\eqref{u-theta}, that, is, a partial regularity result in the spirit of what Caffarelli, Kohn and Nirenberg did for the Navier-Stokes equations in dimension $3$. As in \cite{MR0673830}, this result follows from a simple covering argument and a so-called $\varepsilon$-regularity result.  Exactly as in \cite{MR0673830}, the $\varepsilon$-regularity result relies on the crucial invariance of \eqref{QG2} with respect to the scaling
$$
\theta_r(t,x):= r^{2\alpha -1}\theta(r^{2\alpha}t, rx)\,,\quad u_r(t,x):= r^{2\alpha -1}u(r^{2\alpha}t, rx)\,.
$$ 
Then, the main theorem in \cite{MR4235801} states that for every $L^2$ initial datum, and for every $\alpha\in [\frac{9}{20}, \frac12]$,  there exists an almost everywhere smooth solution of 
\eqref{QG2}-\eqref{u-theta} for which a Hausdorff dimension of the closed set of its singular points is given. More precisely, it is proved that this Hausdorff dimension does not exceed $\frac1{2\alpha}\(\frac{1+\alpha}\alpha(1-2\alpha)+2\)$.

\smallskip
Also as in \cite{MR0673830}, Colombo and Haffter introduced the notion of a well-chosen \it suitable weak solution \rm which allows them to perform energy estimates of nonlinear type controlling a potentially large power of $\theta$. Such nonlinear energy estimates exploit the boundedness of
Leray-Hopf weak solutions in an essential way and are not available for the Navier-Stokes equations.

\smallskip
\section{Jets, cavities and other fluid  interface  problems}\label{sec:jets}

In this section we discuss a number of works of Caffarelli and collaborators where free boundaries arise in the modeling of important physical phenomena, like in the study of the shape and regularity of jets emerging from a nozzle, or rotating fluids in stars or interfaces between fluids flowing together.

\smallskip
Chronologically, the first article was concerned with the study of equilibrium states for simple models of rotating stars or vortex rings \cite{MR0560219}. In both cases the model was a variational problem where the unknown was the density of the rotating fluid. The boundary of the zone occupied by the fluid was a free boundary, and the properties of the fluid region were studied. In the case of rotating stars, the energy functional used to model the rotating ``fluid" in a star was
$$\mathcal E(\rho)=\int_{\R^3}A(\rho(x))dx-\frac12\iint_{\R^3\times\R^3}\frac{\rho(x)\rho(y)}{|x-y|}dx\,dy-\int_{\R^3}J(r)\rho(x)\,dx\,,$$
where $A$ and $J$ were functions related to the material composing the star or the dynamic properties expected for it. 

\smallskip
For physically relevant functions $A$ and $J$, Auchmuty \it et al \rm had proved the existence of a global minimizer of the above energy functional (see \cite{MR0337260, MR0446076}). What Caffarelli and Friedman did in \cite{MR0560219} was to analyse the regularity of the free boundary (of the star shape) and other geometric properties. Among other things they proved that, both in the compressible and the incompressible cases,  the star occupies a compact set made of a finite number of rings (torus-shaped regions) and that it has a smooth boundary. Similar results were proved for the vortex rings model which had been previously studied by Fraenkel and Berger in \cite{MR0422916}.

\smallskip
Then started a very fruitful collaboration of Caffarelli with both  H.W. Alt and A. Friedman  that would last several years and would produce a good number of very interesting works and techniques. There was at the beginning a work on the study of $2$-dimensional jet flows \cite{MR0637494}  and another one  about $2$ and $3$-dimensional infinite cavitational flows past a nose \cite{MR0642623}, followed by a third one on jet flows with gravity \cite{MR0647374}.  In the three cases, a fluid is projected against a hole, or nose, or ``mouth", and what is of interest is the shape of the fluid region past that hole. Early modeling of these phenomena can be found in \cite{MR0088230, MR0047438} and existence of solutions proved in \cite{Leray1935/36, Leray-Weinstein, Weinstein} by Leray and Weinstein,  and in \cite{MR0053684} by Garabedian, Lewy and Schiffer. The main results in these papers are: existence and analysis of the asymptotic direction of the jet without assuming that the nozzle is convex in \cite{MR0637494}, monotonicity and regularity in the case of \cite{MR0642623}. 

\smallskip
In \cite{MR0647374}, the problem of a stationary jet flow for ideal fluids is studied for the first time in the presence of a gravity field. And this is done in the axisymmetric three-dimensional case and in the general two-dimensional one. In both cases it is shown that the jet detaches smoothly at one endpoint of the nozzle.

\smallskip
The authors show that the free boundary problems tackled in \cite{MR0647374} can be stated as variational problems for which they find a solution by approximating them by some cut-off problems. The Maximum Principle, comparison theorems and the Harnack inequality are an important part of the tools used in the proofs. In the end, the uniform estimates for the solutions of the cut-off problems are enough to pass to the limit and recover a solution for the initial problem.

\smallskip
In 1984, Alt, Caffarelli and Friedman published a very important paper about the study of minimization problems with two phases and a free boundary between them \cite{MR0732100}. It was not related to fluids, but was motivated by the study of models for jets and cavities, and it contained a technical result that they would use later in several papers concerning fluid interfaces, and that has had a very big impact in the literature. This was the so-called \it monotonicity lemma\rm. Because it is a very important result that would be used again and again, it will be fully stated below. 

\smallskip
Let $u$ be a minimizer for the energy functional
$$ J[v]:= \int_\Omega \(|\nabla v|^2 + q(x)^2\lambda[v]^2\)dx\,,$$
where $q(x)\neq 0$ for all $x\in \Omega$ and $\lambda[v]$ takes different values for $v$ positive or negative and $\lambda^2[v]$ is lower semicontinuous at $v=0$.

\smallskip
The \it monotonicity lemma \rm states that at any point in the interface, that is, when $u(x_0)=0$, the function
$$r\mapsto r^{-4}\int_{B_r(x_0)} \rho^{2-n}|\nabla u^+|^2dx \cdot  \int_{B_r(x_0)}\rho^{2-n}|\nabla u^-|^2dx $$
is monotone increasing.

\smallskip
This lemma has proved to be a rather powerful tool to establish Lipschitz continuity of the free boundary at well as to identify possible blow-up limits. It has been used very often, in its original form, or in other forms suitable to deal with other  free boundary problems. For instance, in a recent article, \cite{MR4120923}, a new monotonicity lemma is devised for a minimisation problem of the same kind as the above one, but with coefficients, and even singular coefficients, in front of the gradient term in the energy. Of course, this lemma has the same flavor as the old one proved by Alt, Caffarelli and Friedman. The recent article  \cite{MR4142364} contains a long history of the \it monotonicity lemma \rm and its many applications in different contexts. In some sense, this result is the first in a long list of powerful techniques proposed by Caffarelli and collaborators and which will impact strongly the work of a large community later on.

\smallskip
The results in  \cite{MR0732100}, and in particular the \it monotonicity lemma \rm proved in it, were very important tools in three articles written by Alt, Caffarelli and Friedman shortly afterwards: the first one was about the dam problem with two fluids, with two free boundaries, one between one of the fluids and the porous dam wall, and the other one between the fluids themselves. The next ones \cite{MR0733897, MR0740956},  were devoted to the study of jets with two fluids, with again one or two free boundaries depending on the geometry of the problem. In all those articles, the monotonicity lemma, together with application of comparison arguments, played a fundamental role in the proof of the regularity of the interfaces.

\smallskip
Other works belonging to this phase of Caffarelli's career, and to the set of articles devoted to fluid problems, are the study of the abrupt and smooth separation of free boundaries in flow problems \cite{MR0818805}, and another one about the study of compressible flows of jets and cavities \cite{MR0772122}, where Alt, Caffarelli and Friedman solved the problem of existence of axially symmetric compressible subsonic flows of jets and cavities, provided the Mach number is small enough and the flow is inviscid, stationary and irrotational. In the previous works on jets and cavities, the fluids had been considered to be  ``ideal" or incompressible, while here a very different physical situation was considered: adiabatic and isentropic flows. In \cite{MR0772122} a subsonic truncation is performed in order to work with elliptic problems and then, get enough estimates to pass to the limit.

\smallskip
\section{Regularity for the Stefan problem.}\label{sec:Stefan} 

The Stefan problem can be stated in the following way: if $\theta (t, x) $ denotes the temperature of water at a given point $x\in \R^n$ and time $t$, if we introduce a piece of ice inside the water, the ice will start melting, $\theta$ then becoming positive. 

Denote by $\Sigma$ the boundary of the ice, that is, $\Sigma= \partial\{\theta>0\}$. The nonnegative function $\theta$ satisfies 
$$\partial_t\theta - \Delta \theta =0 \quad\mbox{in the water region}\;  \{\theta>0\}\,.$$

The evolution of the free boundary $\Sigma$ will be given by the evolution equation 
 $$\partial_t \theta = |\nabla \theta|^2\;\;\mbox{on}\;\; \Sigma\,.$$ 

In \cite{MR0454350} Caffarelli proved that  $\Sigma$ is infinitely differentiable in $x$ and $t$ outside a closed set of singular points. Moreover,  if one defines $v(x):= \int_0^t \theta(s,x)\,ds$, and if $n=2$, $\theta=v_t$ is continuous across the free boundary $\Sigma$. 

\smallskip
Furthermore, in  \cite{MR0466965} Caffarelli proved that any $C^1$ curve contained inside $\Sigma\times (0, T)$ is necessarily contained in a fixed time slice $\{t=t_0\}$.
The proof was based on a very clever use of the Maximum Principle and the Harnack inequality.

\medskip
Only very recently the above result has been improved: in \cite{MR4695505}, A. Figalli, X. Ros-Oton and J. Serra studied more in detail the singular set for the solution of the one-phase Stefan problem and proved the following:

\smallskip
-- The singular set has parabolic Hausdorff dimension at most $n-1$.

\smallskip
-- The solution admits a $C^\infty$-\,expansion at all singular points, up to a set of parabolic Hausdorff dimension
less than or equal to $n-2$.

\smallskip
-- In $\R^3$, the free boundary is smooth for almost every time $t$, and the set of singular times  has
Hausdorff dimension less than or equal to  $1/2$.

\smallskip
In some later papers, Caffarelli and collaborators treated a more complex ice-water transition model called the $2$-phase Stefan problem, as well as the so-called fractional Stefan problem, where the diffusion is not driven by the linear Laplacian, but by a fractional Laplacian $(-\Delta)^s$ with $s\neq 1$.

\smallskip
\section{Existence and regularity for capillary drop models}\label{sec:drops}

Contact angle hysteresis is an important physical phenomenon, which we can observed in Nature, but that has many important consequences in practice, even in the design of many industrial processes. It is the cause of the strange fact that rain drops stick to window panels instead of sliding down under the action of gravity. The cause of contact angle hysteresis is the chemical and physical heterogeneities of the panels. It is observed and it can be used in practice.

\smallskip
Very interesting results of L.A. Caffarelli and A. Mellet were devoted to the study of liquid drops lying on an inhomogeneous (rough) surface, with or without the action of gravity on the drop. They published the results in \cite{MR2373730, MR2307770}. 

\smallskip
The model they studied is a variational one: for any given set E (the support of the drop), define the following energy functional:
$$\mathcal J(E):= \iint_{z>0}|D\varphi_E|- \int_{z=0} \beta(x,y)\,\varphi_E(x,y,0)\,dx\,dy + \frac1\sigma\int_{z>0}\Gamma\rho\,\varphi_E ,$$
where $(x, y, z)\in \R^3$,  $\varphi_E$ is the characteristic function of $E$, $\sigma$ denotes the surface tension, $\beta$ is the relative adhesion coefficient between the fluid and the solid, $\Gamma$ denotes the gravitational potential and $\rho$ is the local density of the fluid.  The aim is to minimize this energy functional among all sets $E$. The Euler-Lagrange equation for the (local or global) minimizers is 
$$2\bf H = \frac{\Gamma \rho}{\sigma}-\lambda\,,$$
where $\bf H$ denotes the mean curvature of the free surface $\partial E$, and a contact angle condition, known as Young-Laplace's law, which reads:
$$cos\, \gamma =\beta(x, y)\,,$$
where $\gamma$ denotes the angle between the free surface $\partial E$ and the support plane $\{z = 0\}$ along the contact line $\partial(E \cap \{z = 0\})$ (measured within the fluid). The function $\beta$ contains the information about the substrate, and in particular, about the surface inhomogeneities. In these two papers (with $\Gamma=0$ in the first one, $\Gamma>0$ in the second one) the aim is to consider periodic inhomogeneities
$$ \beta = \beta(x/\varepsilon, y/\varepsilon)\,,$$ $\beta$ being $\Z^2$-periodic. Then, study the limit when $\varepsilon\rightarrow 0$. The main result contained in these two articles is that periodic homogenization leads sometimes to a contact angle hysteresis given by 
\be{anglehys}
\cos\, \gamma\in [\gamma_1, \gamma_2], \quad\gamma_1 < \gamma_2\,,
\ee
  at least in the gravity-free case. As explained in \cite{MR2307770}, this means that ``on a perfect surface, the drop should always slide down the plane, no matter how small the inclination. However, a water drop resting on a plane that is slowly inclined will first change shape without sliding, and will only start sliding when the inclination angle reaches a critical value: in the lower parts of the drop, the liquid has a tendency to advance and the contact angle increases until it reaches the advancing contact angle, while in the upper parts of the drop, the liquid has a tendency to recede and the contact angle decreases until it reaches the receding contact angle".

\smallskip
In the first paper published by Caffarelli and Mellet \cite{MR2373730}, global minimizers for the problem without gravity, that is, with $\Gamma=0$, were proven to exist for every $\varepsilon>0$. Moreover, those global minimizers were shown to converge uniformly to a solution satisfying the Young-Laplace's law
$$ \cos \,\gamma = \iint_{[0,1]^2} \beta(x,y)\,dx\,dy\,.$$
So, no contact angle hysteresis in this case, which was actually expected for global minimizers.  But in \cite{MR2307770} some local minimizers are shown to converge towards a solution of the homogenized problem with contact angle hysteresis, that is, satisfying \eqref{anglehys}.

\smallskip
For $\varepsilon>0$ the local minimizers of interest are constructed by solving a constrained minimization problem. The constraining procedure consists in  building a barrier function whenever drops on inclined surfaces are considered. Those barriers prevent the drops to fall and they are so well chosen, that when passing to the limit as $\varepsilon \rightarrow 0$, the limit function is actually an unconstrained local minimizer for the initial problem with contact angle hysteresis. The construction of the barrier functions was quite sophisticated and required a lot of insight.

\smallskip
Note that  Caffarelli had started to work on capillary drops years earlier, see for instance his paper with A. Friedman on the regularity of the boundary of a capillary drop on an inhomogeneous plane and related variational problems \cite{MR0834357}.

\smallskip
\section{ Interpolation inequalities with weights: the Caffarelli-Kohn-Nirenberg inequalities and extensions}\label{sec:CKN}

As already noted in Section \ref{sec:NS}, in Section 7 of \cite{MR0673830} we find a family of inequalities which were used to prove Theorems \ref{thm:THEOREMC} and \ref{thm:THEOREMD} which are the following:
\be{CKN-ineq-orig}
\(\irn{|x|^{-b p} |u|^p}\)^{1/p}\leq C\(\irn{|x|^{-a r}|\nabla u|^r}\)^{c/r}\(\irn{|x|^{\beta q}|u|^q}\)^{(1-c)/q}\,,
\ee
with 
$$ q,r\ge 1\,,\;\;p>0\,,\;\; 0\le c\le 1\,,\;\; -b= c\,\sigma + (1-c)\beta\,,$$

$$\frac1r - \frac{a}{n}\,,\quad \frac1q+\frac\beta{n}\,,\quad \frac1p- \frac{b}{n} >0\,.$$
The main result proved in \cite{MR0673830} and discussed in much more detail in \cite{MR0768824} is the following:

\begin{theorem}[\cite{MR0673830, MR0768824}]\label{thm:CKN-ineq} Let $n\ge 3$. 
The constant $C$ in \eqref{CKN-ineq-orig} is finite if and only if the following relations hold true:

$$
\frac1p-\frac{b}{n}= c\(\frac1r- \frac{a+1}n\)+ (1-c)\(\frac1q+\frac\beta{n}\)\,,
$$

\smallskip
$$ 0\le -a -\sigma\quad \mbox{if}\quad c>0\,, $$
and
$$ -a -\sigma\le 1\quad \mbox{if}\quad c>0 \quad\mbox{and}\quad \frac1r-\frac{a+1}{n}= \frac1p-\frac{b}{n}\,.$$
\end{theorem}
The Caffarelli-Kohn-Nirenberg inequalities (CKN inequalities in short) have attracted a lot of attention in various fields of mathematics, because of their usefulness when studying variational problems and nonlinear PDEs where weights are present as coefficients, and also in the study of blow-up phenomena for evolution nonlinear partial differential equations, in global analysis and in mathematical physics. Also, variants and extensions of these inequalities using $p$-Laplacians, fractional Laplacians and anisotropic weights have been widely studied these last years. 

\smallskip
In the sequel of this section we shall present the main results that have been proved concerning these inequalities, as this is of independent interest, even if not directly concerned with fluids.

\smallskip
Let us define the set of \it admissible functions \rm for \eqref{CKN-ineq-orig} by
$$X:= \{u\in L^q(\R^n, |x|^{\beta q}\,dx) : \irn{|x|^{-a r}|\nabla u|^r}<\infty\}\,.$$

\smallskip
Let us first remark that  the constant $C$ in \eqref{CKN-ineq-orig} satisfies
\be{minpb}
\frac1C = \min_X \,\mathcal F[u]\,,\quad \mathcal F[u]:= \frac{\(\irn{|x|^{-a r}|\nabla u|^r}\)^{a/r}\(\irn{|x|^{\beta q}|u|^q}\)^{(1-a)/q}}{\(\irn{|x|^{-b p} |u|^p}\)^{1/p}}\,,
\ee
and that $v$ is a extremal function for \eqref{CKN-ineq-orig} if and only if it is a global minimizer for \eqref{minpb}.

\smallskip
Moreover, if we define by $X^*$ the space of functions in $X$ that are radially symmetric, and $C^*$ the best constant in \eqref{CKN-ineq-orig} if we restrict the inequality to radially symmetric functions, $v^*$ is a \it radial extremal\rm, that is, an extremal function for \eqref{CKN-ineq-orig} when considering it only among radially symmetric functions if and only if it is a minimizer for the minimization problem
\be{minpbrad}
\frac1{C^*} = \min_{X^*} \, \,\mathcal F[u]\,,
\ee

\smallskip
In \cite{MR1794994} Catrina and Wang studied the case  $r=2, \;c=1$, what without restriction means to consider the set of parameters $\{(a,b): a<\frac{n-2}2, \;a\le b\le a+1\}$. Among other things, they proved the following:

\medskip
-- Whenever $a<b<a+1, \,a<\frac{n-2}2$, there exists an extremal function for \eqref{CKN-ineq-orig} in $X$.

\smallskip
-- When $b=a+1$ or when $a<0, \;b=a$, the constant $C$ is finite, but there is no extremal function for the inequality in $X$.

\smallskip
-- There is a nonpositive number $a_0$ and a continuous and explicit curve $a\mapsto h(a)$ with $h(a_0)=a_0$, $a<h(a)<a+1$ for all $a<a_0$, and $h(a)-(a+1) \rightarrow 0$ as $a\rightarrow -\infty$ such that the extremal functions for \eqref{CKN-ineq-orig} are not radially symmetric if $a<b<h(a)$. 

\smallskip
-- The boundary of the region where the extremals are not radially symmetric contains the set $\{(a,a)\in \R^2: a<0\}$.

\smallskip
-- Up to dilation,  the \it radial extremal \rm is the function
\be{radextremal}
v^*(x):= \( 1+ |x|^{(p-2)(a_c-a)}  \)^{-\frac2{p-2}}\;\;\forall x\in \R^n\,,\quad  \quad a_c:=\frac{n-2}2\,.
\ee

\smallskip
If we say that there is \it symmetry breaking \rm when the extremals for \eqref{CKN-ineq-orig} are not radially symmetric, then Catrina and Wang proved that there is a big subset of the set $\{(a,b)\in \R^2: a<0, \;a<b<a+1\}$ where there is \it symmetry breaking. \rm The method they used to prove it was based on the study of the functional $ \mathcal F[u]$ around the radial extremal function $v^*$: they partially identified the region where $v^*$ is not a local minimum of  $\mathcal F$ in the whole space $X$. That is, there is a big set of parameters $(a,b)$ for which $v^*$ is \it linearly unstable. \rm  In \cite{MR1973285}, Felli and Schneider completed the work of Catrina-Wang, proving that $a_0=0$. 

\smallskip
On the other hand, a series of papers discussed the  \it symmetry set, \rm that is, the sets of  parameters $(a,b)$ for which $v^*$ is a global extremal function, a global minimizer for $\mathcal F$ in $X$.   By putting together results proved in  \cite{MR1223899, MR1731336, MR1734159, MR2560127}, the \it symmetry \rm zone was known to contain the set $\{(a,b)\in \R^2: 0\le a < a_c, \, \;a<b<a+1\}$, the triangle $\{(a,b) : a<0, \,0<b<a+1\}$ and a lower neighborhood of the set $\{(a, a+1): a<0\}$.

\smallskip
Still in the case $r=2, \;c=1$, in \cite{MR3570296}  Dolbeault, Esteban and Loss settled completely the question of \it symmetry \rm and  \it symmetry breaking \rm by proving that \it symmetry breaking \rm takes place if and only if $v^*$ is \it linearly unstable, \rm that is, if and only if $a<0$ and $a<b<h(a)$. 

\smallskip
The above optimal result follows from  proving that for $a<0$, $h(a) \le b<a+1$, there is a unique positive solution to the equation
$$
-\mbox{div}\!\(|x|^{-a}\nabla u\) = |x|^{-bp}|u|^{p-2}u \;\;\mbox{in}\;\;\R^n\,.
$$
This is done by using a flow method, \it \`a la \rm Bakry-Emery, showing that the energy functional 
$$\irn{|x|^{-a r}|\nabla u|^r}$$ is monotone decreasing along a well-chosen highly nonlinear fast-diffusion flow.

\smallskip
The above \it symmetry / symmetry breaking \rm result was expected to hold true. Indeed, it had been conjectured that the only manner to break the symmetry of the global extremal of the inequalities \eqref{CKN-ineq-orig} when $c=1$, $r=2$, would be in cases in which the \it symmetric extremal function \rm would be linearly unstable, making it imposible that it is at the same time a global extremal function. Surprisingly, in \cite{MR2846263}  Dolbeault, Esteban, Tarantello and Tertikas exhibited cases in which the \it symmetric extremal function \rm for \eqref{CKN-ineq-orig} is linearly stable but there is symmetry breaking, that is, there are non radially symmetric solutions 	at which $\mathcal F$ is strictly smaller than at any radially symmetric function. Naturally if this is the case, the exponent $c$ must be strictly less than $1$, because otherwise these examples would contradict the theorem proved in \cite{MR3570296}. The results in \cite{MR2846263} are somehow quantitative in the sense that they provide estimates on the region of parameters for the cases where this anomalous situation arises. Moreover,  \cite{MR3043639} presents a scenario for explaining the reasons why \it  symmetry breaking \rm arises in cases of \it linear stability \rm for the \it radial extremal. \rm

\smallskip
To finish, let us note that many results for the CKN inequalities have been published containing extensions of the results cited above dealing for instance with versions of the inequalities involving fractional and weighted derivatives, anisotropic weights, etc.

\bigskip{\small\noindent\Emph{Acknowledgements.} Special thanks to 
M. Colombo, A.-L. Dalibard, E. Feireisl, I. Gallagher, Ch. Prange and V. Sverak for very useful discussions.
\bigskip\raggedbottom
\bibliographystyle{plain}


\bigskip\begin{center}\rule{2cm}{0.5pt}\end{center}\bigskip

\end{document}